\documentclass{amsart}
\usepackage{amsfonts,amssymb,epsfig}
\usepackage[latin1]{inputenc}
\newtheorem{thm}{Theorem}[section]

\theoremstyle{definition}
\newtheorem{defn}[thm]{Definition}
\theoremstyle{remark}

\numberwithin{equation}{section}

\newcommand{\bbS}{{\mathbb S}}
\newcommand{\bbR}{{\mathbb R}}
\begin{document}

\title[A note on focus-focus singularities, II]{Another note on focus-focus singularities}

\author{Nguyen Tien Zung}
\address{GTA, UMR 5030 CNRS, Département de Mathématiques, Université Montpellier II}
\email{tienzung@math.univ-montp2.fr {\it URL}: www.math.univ-montp2.fr/\~{}tienzung}
\date{third version, Mars/2002, preprint math.DS/0110148, to appear in Lett. Math. Phys.}%

\begin{abstract} {We show a natural relation between the monodromy formula for
focus-focus singularities of integrable Hamiltonian systems and a formula of
Duistermaat-Heckman, and extend the main results of our previous note
\cite{ZungFocus} ($\bbS^1$-action, monodromy, and topological classification) to the
case of degenerate focus-focus singularities. We also consider the non-Hamiltonian
case, local normal forms, etc.}
\end{abstract}
\maketitle

\section{Introduction}

The {\it monodromy} of an integrable Hamiltonian system has been introduced by
Duistermaat in \cite{Duistermaat} as an obstruction to the existence of global
action-angle variables, with a non-trivial example - the spherical pendulum -
provided by Cushman. Since then, this monodromy phenomenon has gained the attention
of many people, from different points of view (see e.g.
\cite{Audin,CuDu,CuSa,San-Focus,WaDu} and references therein). For quite some time,
and even in some very recent papers like \cite{WaDu}, the monodromy had been
calculated mainly by brute force, involving many pages of complicated computations
of action functions. However, most of the integrable systems for which these
computations have been done have a simple common feature : they possess so-called
{\it focus-focus singularities}. And many authors have in fact calculated, case by
case, the monodromy of the Liouville torus fibration around these singularities.
Recall that, a point $x$ on a 4-dimensional symplectic manifold $(M^4,\omega)$ with
an integrable system given by a moment map ${\bf F} = (F_1,F_2) : (M^4,\omega) \to
{\mathbb R}^2$ is called a {\it nondegenerate focus-focus singular point} if
$dF_1(x) = dF_2(x) = 0$, and the quadratic parts of $F_1$ and $F_2$ at $x$ can be
written as $F_1^{(2)} = a(x_1y_1 + x_2y_2) + b (x_1y_2 - x_2y_1), \ F_2^{(2)} =
c(x_1y_1 + x_2y_2) + d(x_1y_2 - x_2y_1)$ in some symplectic system of coordinates
$(x_1,y_1,x_2,y_2)$, with $ad - bc \neq 0$.
%Roughly speaking, when a Hamiltonian
%vector field vanishes at a point, then the spectrum of its linear part is symmetric
%with respect to both the real axis and the pure imaginary axis, and the focus-focus
%case is when we have a quadruple of eigenvalues of the type $(a + ib, a-ib, -a+ib,
%-a -ib)$ with $a,b \in {\mathbb R} \backslash \{0\}$.
The level sets of the moment map form a {\it singular Lagrangian fibration}. If a
singular fiber contains a focus-focus singular point, then one says that it is a
{\it focus-focus singular fiber}. Similarly, one can define corank-2 focus-focus
singularities in integrable Hamiltonian systems with more than two degrees of
freedom (see e.g. \cite{ZungAL}).
%Focus-focus is one of only three types (namely : elliptic, hyperbolic, and
%focus-focus) of components of nondegenerate singularities of integrable Hamiltonian
%systems, see \cite{ZungAL}, and so they are very often seen in practice.
The following simple formula has been obtained in \cite{ZungFocus,ZungAL} : the
monodromy around a nondegenerate focus-focus singular fiber that contains $k \geq 1$
focus-focus singular points (and no hyperbolic singular point) in an integrable
Hamiltonian system with two degrees of freedom is given by the following matrix :
\begin{equation}
\label{eqn:monodromy}
\begin{pmatrix} 1 & k \\ 0 & 1
\end{pmatrix} \ .
\end{equation}

In this note, we will make Formula (\ref{eqn:monodromy}) and other results of our
first note on focus-focus singularities \cite{ZungFocus} easier to use by extending
them to the (possibly) degenerate case. We will also show a direct relation between
Formula (\ref{eqn:monodromy}) and a Duistermaat-Heckman formula.

\begin{defn}
\label{defn:focus} Let $N$ be a compact singular fiber (i.e. singular level set of
the moment map ${\bf F} = (F_1,F_2): (M^4,\omega) \to {\mathbb R}^2$) of an
integrable Hamiltonian system with two degrees of freedom. Then $N$ is called a {\it
possibly degenerate focus-focus singular fiber}, and its singular points are called
{\it possibly degenerate
focus-focus points}, if $N$ satisfies the following  conditions : \\
a) There is a relatively compact saturated connected neighborhood $U(N)$ of $N$ in
which $N$ is the only singular
fiber. \\
b) $N$ contains only a finite number of singular points of the system, and $N$ minus
these singular points is homeomorphic to a non-empty union of cylinders $\bbS^1
\times \bbR^1$.
\end{defn}

{\it Remark.} Of course, if $N$ is a compact singular fiber whose singular points
are all nondegenerate focus-focus, then $N$ will satisfy the conditions of the above
definition, see e.g. \cite{ZungFocus}. Conversely, if $N$ is a possibly degenerate
focus-focus fiber, and $x \in N$ is a nondegenerate singular point, then it is easy
to see that $x$ must be of nondegenerate focus-focus type.

{\it Example.} Consider a Lagrange top under a potential energy field which is not
given by the usual linear function of height, but say by a quadratic function. By
varying the coefficients of this function, one will get degenerate focus-focus
singularities at some parameters.

\begin{thm}
\label{thm:S1} If $N$ is a possibly degenerate focus-focus fiber in an integrable
Hamiltonian system with two degrees of freedom, then there is a neighborhood $V(N)$
of $N$, such that in $V(N)$ there is a Hamiltonian ${\mathbb S}^1$-action with the
following properties : \\
a) This action preserves the system (i.e. the moment map). \\
b) This action is free outside the singular points of the singular fiber $N$, and
fixes these singular points. \\
c) The action has weights $(1,-1)$ at singular points. In other words, near each
singular point of $N$ there is a local symplectic system of coordinates
$(x_1,y_1,x_2,y_2)$ in which the action is generated by the Hamiltonian vector field
$(x_1
\partial/\partial y_1 - y_1
\partial/\partial x_1) - (x_2 \partial/\partial y_2 - y_2 \partial/\partial x_2)$ of
Hamiltonian function $\frac{x_1^2 + y_1^2}{2} - \frac{x_2^2 + y_2^2}{2}$.
\end{thm}

\begin{thm}
\label{thm:monodromy} If $N$ and $U(N)$ satisfy the conditions of Definition
\ref{defn:focus}, then $U(N) \backslash N$ is fibred by regular tori, and the
monodromy of this fibration around $N$ is given by the matrix $\begin{pmatrix} 1 & k
\\ 0 & 1
\end{pmatrix}$, where $k$ is the number of singular points in $N$.
\end{thm}

\begin{thm}
\label{thm:classification} The only topological invariant of a possibly degenerate
focus-focus singularity is its number of singular points. In other words, if $N$ and
$N'$ are two possibly degenerate focus-focus fibers of two integrable Hamiltonian
systems with two degrees of freedom, and the number of singular points in $N$ is
equal to the number of singular points in $N'$, then there is a homeomorphism from a
neighborhood of $N$ to a neighborhood of $N'$ which sends $N$ to $N'$ and which
preserves the singular foliation by Liouville tori.
\end{thm}

The above theorems generalize the main results of \cite{ZungFocus} to the degenerate
case. In particular, Theorem \ref{thm:classification} means that degenerate
focus-focus singularities are homeomorphic to nondegenerate ones. And since
monodromy is a topological invariant, Theorem \ref{thm:monodromy} may be seen as a
consequence of Theorem \ref{thm:classification} and the monodromy formula for
nondegenerate focus-focus singularities. However, we will give in this note a proof
of Theorem \ref{thm:monodromy} which doesn't make use of Theorem
\ref{thm:classification}, but which uses Theorem \ref{thm:S1} a Duistermaat-Heckman
formula with respect to a symplectic $\bbS^1$-action instead. In our opinion, this
natural relation between the monodromy and a Duistermaat-Heckman formula is as
interesting as the monodromy formula (\ref{eqn:monodromy}) itself.

The rest of this note is organized as follows : in Section \ref{section:S1} we prove
Theorem \ref{thm:S1}, and by the way give a new proof of the local normal form, due
to Vey \cite{Vey} and Eliasson \cite{Eliasson}, of nondegenerate focus-focus points.
In Section \ref{section:monodromy} we prove Theorem \ref{thm:monodromy}, and in
Section \ref{section:classification} we prove Theorem \ref{thm:classification}.
Section \ref{section:nonHamiltonian} is devoted to the non-Hamiltonian focus-focus
case, first studied by Cushman and Duistermaat \cite{CuDu}. Various final
observations and remarks are given in Section \ref{section:remarks}, the last
section of this note.

\section{The $\bbS^1$-action}
\label{section:S1}

{\it Proof of Theorem \ref{thm:S1}}. Let us recall how to find the ${\mathbb
S}^1$-action using standard arguments \cite{ZungFocus,ZungAL}. Let $N$ be a
degenerate focus-focus singular fiber, with a saturated connected neighborhood
$U(N)$ as in Definition \ref{defn:focus}. Denote by $F_1$ and $F_2$ the two
components of the moment map, and by $X_1 = X_{F_1}$ and $X_2 = X_{F_2}$ their
corresponding Hamiltonian vector field. Let $m \in N$ be a non-singular point in
$N$. By definition $m$ lies in a cylinder in $N$, implying that there are two real
numbers $k(m),h(m)$ such that the vector field $k(m)X_1 + h(m)X_2$ is periodic of
exact period $2\pi$ at $m$. The numbers $k(m),h(m)$ do not depend on the position of
$m$ inside its cylinder in $N$, though a-priori they may depend on the cylinder
itself (in fact they don't, as will be shown). Let $D$ be a sufficiently small
smooth 2-dimensional disk that intersects $N$ transversally at $m$, and by $V$ the
union of regular Liouville tori which intersect with $D$. Since $U(N)$ is saturated
by assumptions, we have that $V \subset U(N)$. Using the implicit function theorem,
for each point $z$ in $D$ we find a unique pair of numbers $(k(z),h(z))$ close to
$(k(m),h(m))$ such that the orbit starting at $z$ of the vector field $X = k(z)X_1 +
h(z)X_2$ is periodic of exact period $2\pi$ (i.e. the period is exactly $2\pi$ but
not a fraction of it). Consider $k$ and $h$ as functions of $z$. In fact, since the
moment map $(F_1,F_2)$ is regular at $m$, we may take $(F_1,F_2)$ as a coordinate
system on $D$, and so we may consider $k$ and $h$ as functions of two variables
$k(F_1,F_2)$ and $h(F_1,F_2)$ on $D$, and then extend them to a neighborhood of $N$
(which contains $V$) via composition with the moment map. Due to the fact that $X_1$
commutes with $X_2$ (and together generate a standard flat affine structure on each
Liouville torus), the vector field $X = k(F_1,F_2)X_1 + h(F_1,F_2)X_2$ is periodic
of period exactly $2\pi$ not only at the points in $D$, but in the whole $V$. In
order to show that $X = k(F_1,F_2)X_1 + h(F_1,F_2)X_2$ is periodic in a neighborhood
of $N$, it suffices to show that the closure of $V$ contains an open neighborhood of
$N$.

Let $y$ be a point which is sufficiently close to $N$ but which does not belong to
$N$. Since the codimension of $N$ in $U(N)$ is $2$, we have that $U(N) \backslash N$
is connected, and therefore there is a path $\gamma$ going from $y$ to a point $z$
inside $D$ ($z$ close to $m$), such that $\gamma(0) = y, \gamma(1) = z$, and for any
$t \in [0,1]$ we have that $\gamma(t)$ is close to $N$ but does not belong to $N$.
Then one checks that the whole path $\gamma$ lies in $V$, because the intersection
of $\gamma$ with $V$ is an non-empty open closed subset of $\gamma$ : Open because
if $\gamma(t) \in V$ for some $t$, then the regular Liouville torus which contains
$\gamma(t)$ will intersect $D$ at a point inside $D$ (because the value of the
moment map at $\gamma(t)$ is very close to its value at $N$, i.e. at $m$), and hence
nearby Liouville tori (which form an open neighborhood of this particular Liouville
torus) also belong to $V$. Closed because if $\gamma(t)$ lies in the closure of $V$,
then there are Liouville tori in $V$ which lie arbitrarily close to the regular
Liouville torus which contains $\gamma(t)$ (by regularity, if a point lies near
$\gamma(t)$ then the Liouville torus which contains this point is everywhere close
to the Liouville torus which contains  $\gamma(t)$). It implies that there are
points in $D$ which are arbitrarily close to the Liouville torus which contains
$\gamma(t)$. In particular, there is a point in the closure of $D$ which lies in the
same Liouville torus as $\gamma(t)$. But considering the fact that the value of the
moment map at this point, which is the same as its value at $\gamma(t)$, is very
near its value at $m$, this point must in fact lie inside $D$, which means that
$\gamma(t)$ belongs to $V$. Thus the whole path $\gamma$ belongs to $V$, and in
particular $y = \gamma(0) \in V$, for any $y$ close enough to $N$ but not belonging
to $N$. We have proved that the closure of $V$ contains a neighborhood of $N$.

The vector field $X = k(F_1,F_2)X_1 + h(F_1,F_2)X_2$ has been shown to be periodic
in a neighborhood of $N$, and is of period exactly $2\pi$ outside $N$. It is clear
that this vector field fixes the singular points of $N$, where $X_1 = X_2 = 0$. Let
us show that at regular points in $N$, the period of $X = k(F_1,F_2)X_1 +
h(F_1,F_2)X_2$ is also exactly $2\pi$. Let $m'$ be a regular point of $N$. Then by
construction, the value of functions $k$ and $h$ at $m'$ is the same as their value
at $m$. If the period of $X$ at $m'$ is smaller than $2\pi$, then it means that
there is a positive integer $p  \geq 2$ such that $k(m)X_1 + h(m)X_2$ is periodic of
period $2\pi/p$ at $m'$. Repeat the above process of finding $X$, but starting at
point $m'$ instead of at point $m$, and with period $2\pi/p$ instead of $2\pi$. In
the end we will get that $k(m)X_1 + h(m)X_2$ is periodic of period $2\pi/p$ (or
smaller) at $m$, which is contradictory to the assumptions about $k(m)$ and $h(m)$.
Thus the period of $X$ at $m'$ must be exactly $2\pi$ also.

By integrating the above periodic vector field $X$, we obtain an $\bbS^1$-action in
a neighborhood of $N$. The fact that this vector field is Hamiltonian follows from
(the proof of) Arnold-Liouville theorem on the existence of action-angle variables,
see \cite{ZungAL}. It is clear that this ${\mathbb S}^1$-action preserves the
system, is smooth (resp., analytic, finitely differentiable) if the system is smooth
(resp., analytic, finitely differentiable), and it is free outside the singular
points of $N$.

Near each singular point $x \in N$, the above ${\mathbb S }^1$-action can be
linearized, i.e. there is a symplectic system of coordinates $(x_1,y_1,x_2,y_2)$ in
which the action is generated by the Hamiltonian function $f = \frac{a(x_1^2 +
y_1^2)}{2} + \frac{b(x_2^2 + y_2^2)}{2}$, i.e.
$$
X = a(x_1 \partial/\partial y_1 - y_1
\partial/\partial x_1) +  b(x_2 \partial/\partial y_2 - y_2 \partial/\partial x_2),
$$
where $a,b \in {\mathbb Z}$. If say $|a| \geq 2$ then the ${\mathbb S}^1$-action is
not free outside $x$ (it has points of period $2\pi/|a|$), and if $a = 0$ then $x$
is not an isolated fixed point of the action. The remaining cases are $a,b = \pm 1$.
If $(a,b) = (1,1)$ or $(a,b) = (-1,-1)$ (which is the same up to an orientation)
then the Hamiltonian $f$ generating the ${\mathbb S}^1$-action has compact level
sets near $x$ (the $3$-spheres), and it's easy to check that on each of these
$3$-spheres the must exist singular points of the integrable system, due to the fact
a $3$-sphere cannot be foliated by regular tori. Thus the only
possibility is $(a,b) = (1,-1)$, up to a permutation. $\diamondsuit$ \\

{\it Remark.} In fact, most examples provided in \cite{ZungFocus} of systems with
focus-focus singularities admit a global ${\mathbb S}^1$-action, except perhaps the
Manakov system on $so(4)$ (a free 4-dimensional rigid body). \\

As a corollary of the existence of an $\bbS^1$-action, we propose here to give a new
proof of the following local normal form theorem due to Eliasson\footnote{See
Eliasson's Ph.D. thesis dated 1984, rather than his paper dated 1990, for his proof.
In the analytic case this result is due to Vey \cite{Vey}.} \cite{Eliasson}.

\begin{thm}[\cite{Eliasson}]
\label{thm:normalform} Let $m \in (M^4,\omega)$ be a focus-focus singular point in a
smooth  integrable Hamiltonian system with two degree of freedom given by a moment
map ${\bf F} = (F_1,F_2): (M^4,\omega) \to {\mathbb R}^2$. Then there exists a local
smooth  symplectic system of coordinates $(x_1,y_1,x_2,y_2)$ near $m$ (the
symplectic form is $\omega = dx_1 \wedge dy_1 + dx_2 \wedge dy_2$) in which the
local singular Lagrangian fibration of the system is given by the two quadratic
functions $f_1 = x_1y_1 + x_2y_2$ and $f_2 = x_1y_2 - x_2y_1$ (in other words, we
have $dF_1 \wedge df_1 \wedge df_2 = dF_2 \wedge df_1 \wedge df_2 = 0$).
\end{thm}

{\it Proof}. It is a corollary of the following three facts: a) existence of a local
Hamiltonian ${\mathbb S}^1$-action \cite{ZungFocus,ZungAL}; b) existence of a formal
Birkhoff normalization (the classical result of Birkhoff); and c) the Hamiltonian
equivariant Sternberg theorem (for an ${\mathbb S}^1$-invariant vector field)
\cite{BeKo,Chaperon}.

Using the existence of a formal Birkhoff normalization, and Borel's theorem on the
approximation of formal mappings by smooth mappings, we can normalize the system up
to a flat term. In other words, we may assume that
$$
F_1 = G_1(f_1,f_2) + \ {\rm flat}, \ F_2 = G_2(f_1,f_2) + \ {\rm flat}
$$
($\rm flat$ means a flat term), where $(G_1,G_2) : {\mathbb R}^2 \to {\mathbb R}^2$
is a local smooth map, and
$$
f_1 = x_1y_1 + x_2y_2, \ f_2 = x_1y_2 - x_2y_1
$$
is the standard quadratic focus-focus moment map in a smooth symplectic system of
coordinates $(x_1,y_1,x_2,y_2)$. The symplectic form is $\omega_0 = dx_1 \wedge dy_1
+ dx_2 \wedge dy_2$.

By the nondegeneracy assumptions, this map is a local diffeomorphism. Since we are
interested in the singular Lagrangian fibration but not the moment map itself, we
can compose $(F_1,F_2)$ with the inverse of $(G_1,G_2)$ to get a new moment map
which satisfies
$$
F_1 = f_1 + \ {\rm flat}, \ F_2 = f_2 + \ {\rm flat} .
$$
Notice that the Hamiltonian vector field of $f_2 = x_1y_2 - x_2y_1$ is periodic. In
this normal form, the Hamiltonian ${\mathbb S}^1$-action around the focus-focus
singular point will be also generated by a function of the type $f_2 + \ {\rm
flat}$, so we may assume (by changing the moment map without changing the fibration)
that $F_2$ is in fact the generator of our ${\mathbb S}^1$-action. Since $F_2 = f_2
+ \ {\rm flat}$, this action is formally equal to the linear action generated by
$f_2$. Thus, using averaging, we can linearize this action by a symplectomorphism
which is formally equal to the identity map. Using this linearization, we can assume
that:
$$
F_1 = f_1 + \ {\rm flat}, \ F_2 = f_2 .
$$
The Hamiltonian vector field $X_1$ of $F_1$ is now formally equal to the linear
Hamiltonian vector field $X_{f_1}$, it is invariant under the linear ${\mathbb
S}^1$-action degenerated by $X_2 = (x_1 \partial/\partial x_2 - x_2
\partial/\partial x_1) + (y_1 \partial/\partial y_2 - y_2 \partial/\partial y_1)$,
and it does not have purely imaginary eigenvalues. Hence one can apply the
equivariant Hamiltonian Sternberg theorem with respect to an ${\mathbb S}^1$-action
\cite{BeKo,Chaperon} to find a local smooth symplectomorphism that maps $X_1$ to
$X_{f_1}$ (and $X_2$ to itself, of course). Thus under this symplectomorphism, the
system becomes linear:
$$
\begin{array}{l}
X_1 = (x_1 \partial/\partial x_1 + x_2
\partial/\partial x_2) - (y_1 \partial/\partial y_1 + y_2 \partial/\partial y_2), \\
X_2 = (x_1 \partial/\partial x_2 - x_2
\partial/\partial x_1) + (y_1 \partial/\partial y_2 - y_2 \partial/\partial y_1) ,
\end{array}
$$
and the singular Lagrangian foliation is given by two quadratic functions $f_1 =
x_1y_1 + x_2y_2, \ f_2 = x_1y_2 - x_2y_1.$ $\diamondsuit$
\\

\section{Monodromy and a Duistermaat-Heckman formula}
\label{section:monodromy}

{\it Proof of Theorem \ref{thm:monodromy}}. We will use the notations introduced in
the previous sections : $N$, $U(N)$, ... The proof will use the Duistermaat-Heckman
measure \cite{DuHe} to show the role played by the ${\mathbb S}^1$-action given by
Theorem \ref{thm:S1}. Without loss of generality, we can assume that $F_2$ (the
second component of the moment map) is the generator of this ${\mathbb S}^1$-action.
For simplicity, we will localize the degenerate focus-focus singularity. Using
integrable surgery \cite{ZungIntegrable2}, we can add elliptic singularities to the
boundary of $U(N)$ (provided $U(N)$ is chosen appropriately) to get an integrable
system on a compact symplectic 4-manifold with the following property : on each
compact level set $Q_c = \{ F_2 = c \}$ with $c \neq 0$ and near $0$, the ${\mathbb
S}^1$-action generated by the periodic vector field $X_2 = X_{F_2}$ is free, and the
set of singular points of the system on $Q_c$ consists of two circles $S_{1,c},
S_{2,c}$ of nondegenerate corank-1 elliptic singularities. On the singular level set
$Q_0 \supset N$, the ${\mathbb S}^1$-action is still free outside the set of
singular points in $N$, and the set of singular points of the system on $Q_0$
consists of two circles $S_{1,0}, S_{2,0}$ of nondegenerate corank-1 elliptic
singularities plus the singular set of $N$. Additionally, we can assume that, for
some small number $\delta > 0$, $\delta$ is the maximum value of $F_2$ and $-\delta$
is the minimum value of $F_2$, and the sets $Q_{\delta}, Q_{-\delta}$ are
two-dimensional sphere consisting of elliptic singular points of the system. For the
${\mathbb S}^1$-action, there are $k$ fixed points of weights $(1,-1)$ (lying in
$N$), two spheres of fixed points $Q_{\delta}$ and $Q_{-\delta}$, and the action is
free outside these sets. The orbit space of the localized integrable system is
topologically a square with a particular (degenerate focus-focus) singular point
inside. The edges of the square correspond to corank-1 nondegenerate elliptic
singularities, and the corners of the square correspond to nondegenerate corank-2
elliptic singularities.

Denote by $V(c)$ the symplectic volume (area) of the quotient of $Q_c$ by the
${\mathbb S}^1$-action, with the reduced symplectic form. Then a special case of a
formula due to Duistermaat and Heckman \cite{DuHe} measures how the behavior of the
function $V(c)$ changes when one passes though the fixed points of the ${\mathbb
S}^1$-action (here the point $0$ in the variable $c$):

%??????  <<--- add details here !!!

$V(c)$ is piecewise linear (linear when $\delta \geq c \geq 0 $ and when $ - \delta
\leq c \leq 0$), and
\begin{equation}
\label{eqn:DH} V(c) + V(-c) = 2V(0) - kc .
\end{equation}
(The minus sign on the right hand side is due to the fact that the weights at the
fixed points are $(1,-1)$). On the other hand, the orbit space of our localized
integrable system has a stratified integral affine structure \cite{ZungIntegrable2},
and $V(c)$ can be viewed as the affine length of the interval $\{ F_2 = c \}$ on the
orbit space with respect to this integral affine structure. Formula \ref{eqn:DH}
then explains the change in the behavior of the integral affine structure in the
orbit space when one passes through the focus-focus point: around this point, the
integral affine structure can be obtained from the standard flat structure in
${\mathbb R}^2 = \{(x, y)\}$ near the origin $O$ by cutting out the angle
$\angle\{(0,1),(-k,1)\}$ and gluing the edges of the rest together by the integral
linear transformation $(x, y) \mapsto (x + ky, y)$. And this change corresponds to
the monodromy formula (\ref{eqn:monodromy}) in action variables, which by duality is
the same thing as the monodromy in the torus
fibration (see e.g. \cite{Duistermaat,ZungIntegrable2}). $\diamondsuit$ \\

\section{Topological classification}
\label{section:classification}

{\it Proof of Theorem \ref{thm:classification}}. Denote by $P$ and $V(P)$ the
quotient of $N$ and $U(N)$ by the Hamiltonian $\bbS^1$-action given by Theorem
\ref{thm:S1}. Due to the fact that the action has weights $(1,-1)$ at each fixed
point, $V(P)$ is a topological 3-manifold, and $P$ is a simple closed curve in it.
Since the moment map is preserved by the $\bbS^1$-action, it can be projected to a
map $(\hat{F_1},\hat{F_2}) : V(P) \to \bbR^2$, which turns $V(P)$ into a
topologically trivial circle bundle. (The singular points in $P$ are singular in the
sooth sense, but non-singular in the homeomorphic sense). Thus if $(N,U(N))$ and
$(N',U'(N'))$ denote two possibly degenerate focus focus singularities with the same
number of singular points (which we will denote by $k$), then at least their
reductions by the $\bbS^1$-action are topologically equivalent : there is a
homeomorphism $\psi$ from $V(P)$ to $V'(P')$ which preserves the circle fibrations,
and which sends the singular points in $P$ (i.e. the images of the singular points
in $N$ under projection) to the singular points in $P'$.

Now one can lift the map $\psi : V(P)\to V'(P')$  to an $\bbS^1$-equivariant map
from $U(N)$ to $U'(N')$: first do it for $N$, then for small neighborhoods of
singular points of $N$, then extend to the rest of $U(N)$. It can be done by using
``local sections'', and is an exercise in elementary topology. There are no
obstructions. $\diamondsuit$

{\it Remark.} In \cite{ZungFocus}, I wrongly suggested that the above homeomorphism
in the nondegenerate case can always be made smooth, i.e. the number $k$ of
focus-focus singular points in the focus-focus singular fiber is the only smooth
invariant - after a more careful analysis, Bolsinov found out that there are some
other smooth invariants in the case $k > 1$ (see \cite{BoFo-Book}). A topological
classification of more general nondegenerate singularities of integrable Hamiltonian
systems is given in \cite{ZungAL}.

\section{The non-Hamiltonian case}
\label{section:nonHamiltonian}

In \cite{CuDu}, Cushman and Duistermaat generalized the results of \cite{ZungFocus}
to the case of  {\it integrable non-Hamiltonian systems with focus-focus
singularities}. They assign to each focus-focus point a sign, either plus or minus,
depending on some orientation (in the Hamiltonian case, all focus-focus points have
plus sign). They again obtained the existence of an ${\mathbb S}^1$-action (which is
not surprising), and the same monodromy matrix as given by Formula
\ref{eqn:monodromy}, but with $k$ now being the number of positive focus-focus
points minus the number of negative focus-focus points. Naturally, the results in
the previous sections of this note can also be extended to the non-Hamiltonian case,
so we also get a generalization of Cushman-Duistermaat's results.

\begin{defn}
\label{defn:focus_nh} If we have two vector fields $X_1,X_2$ and two functions
$F_1,F_2$ on a 4-manifold $M^4$ such that $[X_1,X_2] = 0, X_1(F_1) = X_1(F_2) =
X_2(F_1) = X_2 (F_2) = 0$, then we say that we have an integrable non-Hamiltonian
system of bi-index $(2,2)$ (i.e. 2 commuting fields and 2 common first integrals). A
point $x \in M^4$ is called singular for such an integrable system, if $X_1 \wedge
X_2 (x) = 0$ or $dF_1 \wedge dF_2 (x) = 0$. A connected common level set $N$ of the
first integrals $F_1,F_2$ is called a {\it possibly degenerate focus-focus singular
fiber}, and its singular points called {\it possibly degenerate focus-focus singular
points}, if $N$ satisfies the following conditions : \\
a) There is a relatively compact neighborhood $U(N)$ of $N$, which is saturated by
connected common level sets of $F_1$ and $F_2$.\\
b) $U(N)$ contains only a finite number of singular points of the system and they
all lie in $N$, and $N$ minus these singular points is homeomorphic to a non-empty
union of cylinders $\bbS^1 \times \bbR^1$. \\
\end{defn}

\begin{thm}
\label{thm:S1_nh} If $N$ is a possibly degenerate focus-focus fiber in an integrable
non-Hamiltonian system of bi-index $(2,2)$, then in a neighborhood $U(N)$ of $N$
there is an ${\mathbb S}^1$-action with the following properties : \\
a) This action preserves the system (i.e. the two commuting vector fields and the
two first integrals) \\
b) This action is free outside the singular points of the singular fiber $N$, and
fixes these singular points. \\
c) The action has weights $(1,\pm 1)$ at singular points. In other words, near each
singular point of $N$ there is a local symplectic system of coordinates
$(x_1,x_2,x_3,x_4)$ in which the action is generated by the  vector field $(x_1
\partial/\partial x_2 - x_2
\partial/\partial x_1) \pm (x_3 \partial/\partial x_4 - x_4 \partial/\partial x_3)$.
\end{thm}

The proof of Theorem \ref{thm:S1_nh} is absolutely similar to that of Theorem
\ref{thm:S1}. The reason why we write weights $(1,\pm 1)$ instead of $(1,-1)$ in the
above theorem is that there is no distinction between $(1,1)$ and $(1,-1)$, before
an orientation is introduced. Note that $U(N)$ is orientable. For example, one can
choose an orientation as follows: Choose a Riemannian metric in $U(N)$. If $m \in
U(N)$ is regular point of the system, then the vectors $X_1(m),X_2(m),\nabla
F_1(m),\nabla F_2(m)$ (where $\nabla$ denotes the gradient of a function) are
linearly independent at $m$, and one can use these vector fields to orientate
$U(N)$.

\begin{defn} Suppose that $N$ is a possibly degenerate focus-focus singular fiber as in
Definition \ref{defn:focus_nh}, and that the ambient 4-manifold $M^4$ is oriented.
Then a singular point $x \in N$ is called {\it positive} if there is a positively
oriented local system of coordinates $(x_1,x_2,x_3,x_4)$ in which the
$\bbS^1$-action given by Theorem \ref{thm:S1_nh} is generated by $(x_1
\partial/\partial x_2 - x_2
\partial/\partial x_1) - (x_3 \partial/\partial x_4 - x_4 \partial/\partial x_3)$,
and {\it negative} otherwise.
\end{defn}

{\it Remark.} Of course, a singular point in $N$ is either positive or negative but
never both. The above definition is an extension of a definition of Cushman and
Duistermaat \cite{CuDu} to the case of possibly degenerate focus-focus points. It is
also clear that in the Hamiltonian case, all  possibly degenerate focus-focus
singular points are of plus sign. If we change the orientation of the ambient
manifold $M^4$, then positive points become negative and vice versa.

\begin{thm}
\label{thm:monodromy_nh} If $N$ and $U(N)$ satisfy the conditions Definition
\ref{defn:focus_nh}, then $U(N) \backslash N$ is fibred by regular tori, and the
monodromy of this fibration around $N$ is given by the matrix $\begin{pmatrix} 1 & k
\\ 0 & 1 \end{pmatrix}$, where $k$ is the number of positive singular points in $N$ minus
the number negative singular points in $N$.
\end{thm}

Theorem \ref{thm:monodromy_nh} may be viewed as a generalization of Theorem
\ref{thm:monodromy} to the non-Hamiltonian case, or also as a generalization of the
results of \cite{CuDu} to the degenerate case. It can be proved, for example, by
combining the results of \cite{CuDu} and the non-Hamiltonian version of Theorem
\ref{thm:classification} : In the non-Hamiltonian case, the full topological
invariant is made up of the cyclic order of the signs of the singular points in $N$.

\section{Some final remarks}
\label{section:remarks}

\subsection{Topological torus fibrations}

From the topological point of view, the above topological classification and
monodromy formulas (both in the Hamiltonian case and the non-Hamiltonian case) may
be seen, after some preparatory results, as special cases of the theory of singular
torus fibrations in 4-manifolds developed by Matsumoto and other people, see
\cite{Matsumoto} and references therein. This theory is in turn inspired by
Kodaira's theory of elliptic fibrations in complex surfaces. Unfortunately or
fortunately, I didn't know of Matsumoto's paper when writing \cite{ZungFocus} -
otherwise that note would have been written differently, or would not appear at all.
Focus-focus singular points correspond to what are denoted in \cite{Matsumoto} as
$I_+$ (the case of positive sign) and $I_-$ (the case of negative sign, which can
only happen in non-Hamiltonian systems). The ${\mathbb S}^1$-action does not seem to
be present in the above-mentioned work of Matsumoto nor in the earlier papers by
Kodaira; they didn't need it. However, I would like to stress here again that this
action (and local torus actions in general) is very useful in the study of
integrable systems. In \cite{ZungFocus}, we showed how to use this action to perturb
an integrable system with two degrees of freedom into an integrable system whose
focus-focus singular fibers contain exactly one singular focus-focus point. Of
course, the same result holds in the degenerate focus-focus case.

\subsection{Monodromy without focus-focus singularities}

Locally, among nondegenerate singularities, only focus-focus components create
non-trivial monodromy (elliptic and hyperbolic components don't). But globally,
there may be domains in the orbit space of a nondegenerate integrable Hamiltonian
system, which does not border any focus-focus singular point, and which still has
non-trivial monodromy. A simple example is the following : consider a spherical
pendulum of radius 1, centered at the origin of $\bbR^3$, not with the standard
potential energy function $z$ (where $z$ is the coordinate of the vertical axis),
but a potential energy function of the type $z - z^2/2R$, where $R$ is a positive
constant sightly smaller than 1. The maximum of the potential is achieved at the
level $z = R$, slightly below the highest possible level ($z=1$) of the pendulum.
Such a spherical pendulum doesn't have a focus-focus point, but still has monodromy
similar to the usual spherical pendulum.

\subsection{Systems with many degrees of freedom}

It follows easily from the results of \cite{ZungAL} that, if $N$ is a corank-2
focus-focus singular fiber of an integrable Hamiltonian system with $n$ degrees of
freedom, $n > 2$, with some additional nondegeneracy condition (called ``topological
stability'' in \cite{ZungAL}) then a tubular neighborhood of $N$ together with the
singular Lagrangian fibration is homeomorphic to the direct product of a tubular
neighborhood of a regular torus in an integrable system with $n-2$ degrees of
freedom with a tubular neighborhood of a focus-focus singular fiber of an integrable
system with two degrees of freedom. Thus, corank-2 focus-focus singularities in
higher dimensions are the same as focus-focus singularities in dimension $4$. for
example, in the case with 3 degrees of freedom, the monodromy matrix is
\begin{equation}
\begin{pmatrix} 1 & k & 0 \\ 0 & 1 & 0 \\ 0 & 0 & 1
\end{pmatrix} \ .
\end{equation}
Still, one must be careful about the bases in which these matrices are written. For
example, Bates and Zou in a paper published in 1993 made a mistake of using $2
\times 2$ matrices when studying monodromy around two focus-focus singularities in a
3-degrees-of-freedom system, and arrived at a wrong multiplication formula.

\subsection{Integral affine structure of the base space}

By duality, the monodromy of the torus fibration of an integrable Hamiltonian system
is the same as the monodromy of the integral affine structure on the base space. In
fact, when we consider the whole base space with its singularities, then it has a
stratified integral affine structure, and a {\it affine monodromy sheaf} associated
with it, see \cite{ZungIntegrable2}. This sheaf is a free Abelian sheaf which is
locally constant on each stratum, and the structure of this sheaf is a topological
invariant of the system. One can argue that, by Bohr-Sommerfeld rules,  the quantum
joint spectrum of a semi-classical integrable system is nothing but a discretization
of the integral affine structure of the base space of the corresponding classical
integrable system, and that's the main reason why quantum integrable systems have
the same monodromy as its classical counterpart, as observed in many places (see
e.g. \cite{CuSa,San-Focus,WaDu} and references therein). In fact, the quantum joint
spectrum mimics the integral affine structure not only near focus-focus
singularities, but near elliptic and hyperbolic singularities as well.

\subsection{Focus-focus singularities in geometry}

The Strominger-Yau-Zaslow conjecture in mirror symmetry \cite{StYaZa} says that on
each ``dualizable'' Calabi-Yau manifold there is a special Lagrangian fibration and
the dual Calabi-Yau manifold corresponds to the dual special Lagrangian fibration.
From the point of view of integrable Hamiltonian systems, these fibrations are
integrable systems with corank-2 focus-focus singularities, and higher-corank
singularities which must be similar in some sense to focus-focus ones. Hence, a good
understanding of focus-focus singularities and its higher-dimensional degenerate
sisters may be helpful in the SYZ construction of mirror symmetry. Though many
authors have written many papers about the SYZ construction, at present the
situation is still not clear, especially in what concerns singularities, to my very
limited knowledge.

Another application of focus-focus singularities is in the geometry of 4-manifolds.
For example, using a symplectic (integrable) surgery involving focus-focus
singularities, Symington constructed so-called {\it generalized symplectic rational
blow-downs}, which are useful for creating interesting 4-manifolds and calculating
their invariants, see \cite{Symington} and references therein.

\vspace{0.5cm} {\bf Acknowledgements}. I would like to thank Richard Cushman and
Michèle Audin for inciting me to write this note and for various discussions on
integrable systems related to this work, and Pierre Molino for stimulating
conversations on smooth normal form problems. I'm also grateful to Hans Duistermaat,
Vu Ngoc San, and the anonymous referee for their critical remarks on earlier
versions of this note.

\bibliographystyle{amsplain}

\end{document}